\documentclass{article}

\usepackage{amsmath,amsfonts,amssymb,a4,enumerate,epsfig,theorem,psfrag}
\usepackage[latin1]{inputenc}

\newcommand{\trace}{{\operatorname{tr}}}
\newcommand{\transpose}{{\ast}}

\newcommand{\NN}{{\mathbb{N}}}

\newcommand{\RR}{{\mathbb{R}}}

\newcommand{\ww}{{\mathbf{W}}}
\newcommand{\ttt}{{\mathbf{t}}}
\newcommand{\Aa}{{\cal{A}}}

\newcommand{\Dd}{{\cal{D}}}
\newcommand{\Ee}{{\cal{E}}}

\newcommand{\Jj}{{\cal{J}}}
\newcommand{\Kk}{{\cal{K}}}

\newcommand{\Ss}{{\cal{S}}}
\newcommand{\Tt}{{\cal{T}}}
\newcommand{\Uu}{{\cal{U}}}
\newcommand{\Vv}{{\cal{V}}}
\newcommand{\Xx}{{\cal{X}}}
\newcommand{\Yy}{{\cal{Y}}}
\newcommand{\Zz}{{\cal{Z}}}
\newcommand{\bT}{{\mathbf{T}}}

\newcommand{\ILa}{\Tt_\Lambda}
\newcommand{\ILai}{\Tt_{\Lambda_i}}
\newcommand{\ILaie}{\Tt_{\Lambda_i} \times (-\epsilon_Z, \epsilon_Z)}
\newcommand{\ILaiei}{\Tt_{\Lambda_i} \times (-\epsilon', \epsilon')}
\newcommand{\JLa}{\bar\Jj_\Lambda}

\newcommand{\UpLa}{\Uu^\pi_\Lambda}
\newcommand{\bQ}{{\mathbf{Q}}}
\newcommand{\bR}{{\mathbf{R}}}

\newcommand{\proof}{{\noindent\bf Proof: }}

\def\qed{\unskip\nobreak\hfil\penalty50\hskip1.75em\null\nobreak\hfil
$\blacksquare$ {\parfillskip=0pt \finalhyphendemerits=0 \par}\medbreak}

\newcommand\capsize{\relax}

\newcommand\diag{\operatorname{diag}}

\newtheorem{lemma}{Lemma}[section]
\newtheorem{theo}[lemma]{Theorem}
\newtheorem{prop}[lemma]{Proposition}
\newtheorem{coro}[lemma]{Corollary}

\title{The asymptotics of Wilkinson's shift iteration}
\author{Ricardo S. Leite, Nicolau C. Saldanha and Carlos Tomei }

\begin{document}
\maketitle

\begin{abstract}
We study the rate of convergence of Wilkinson's shift iteration
acting on Jacobi matrices with simple spectrum. We show that for
{\it AP-free} spectra (i.e., simple spectra containing no arithmetic
progression with $3$ terms), convergence is cubic. In order $3$,
there exists a tridiagonal symmetric matrix $P_0$ which is the limit
of a sequence of a Wilkinson iteration, with the additional property
that all iterations converging to $P_0$ are strictly quadratic.
Among tridiagonal matrices near $P_0$, the set $\Xx$ of initial
conditions with convergence to $P_0$ is rather thin: it is a union
of disjoint arcs $\Xx_s$ meeting at $P_0$, where $s$ ranges over the
Cantor set of sign sequences $s: \NN \to \{1,-1\}$. Wilkinson's step
takes $\Xx_s$ to $\Xx_{s'}$, where $s'$ is the left shift of $s$.
Among tridiagonal matrices conjugate to $P_0$, initial conditions
near $P_0$ but not in $\Xx$ converge at a cubic rate.
\end{abstract}

\medbreak


\smallbreak

{\noindent\bf MSC-class:} 65F15; 37E30.

\section{Introduction}

In this paper, we study of the asymptotics of the shifted $QR$
iteration with the so called \textit{ Wilkinson's shift}, acting on
Jacobi matrices. More precisely (\cite{Wilkinson}, \cite{Demmel},
\cite{Parlett}), for an $n \times n$ real symmetric tridiagonal
matrix $T$ and a real number $s$, write, if possible, the unique
{\it $QR$ factorization} $T-sI = Q R $, where $Q$ is orthogonal and
$R$ is upper triangular with positive diagonal. Wilkinson's shift
strategy is the choice of $s = \omega(T)$ equal to the eigenvalue of
the bottom $2 \times 2$ principal minor of $T$ which is closest to
the bottom entry $T_{nn}$. A \textit{Wilkinson step} obtains a new
matrix
\[\ww(T) =  Q^\transpose T Q = R T R^{-1}. \]
From both defining formulae, $\ww(T)$ is symmetric and upper
Hessenberg, and thus, must also be a real, symmetric tridiagonal
matrix, with the same spectrum as $T$, as well as the signs of the
nontrivial off-diagonal elements.

As is well known (\cite{Parlett}), if the iterates $T_k = \ww^k(T)$
of the Wilkinson step starting from a Jacobi matrix with simple
spectrum exist for arbitrary $k$, then their lowest off-diagonal
entry tend to $0$: we are interested in the rate of convergence of
this sequence. It has been conjectured (\cite{Parlett},
\cite{Demmel}) that the rate is cubic, i.e., for any Jacobi matrix
$T$ there exists a constant $C$ such that $|(T_{k+1})_{n,n-1}| \le C
|(T_{k})_{n,n-1}|^3$. As we shall see, this is true for most
matrices $T$ but false in general. A Jacobi matrix $T$ is an {\it
AP-matrix} if its spectrum contains an arithmetic progression with
three terms and is {\it AP-free} otherwise. For AP-free matrices,
the conjecture is indeed true. On the other hand, there exist $3
\times 3$ AP-matrices for which the rate of convergence is, in the
words of Parlett, merely quadratic.

We give an outline of the proof.
A basic ingredient are the {\it bidiagonal
coordinates}, consisting of eigenvalues  $\lambda_i$, $i = 1,
\ldots, n$, and additional variables $\beta^\pi_i$, $i = 1, \ldots,
n-1$, defined on large open sets of tridiagonal matrices. As proved
in \cite{LST}, the set $\ILa$ of real, symmetric, tridiagonal
matrices with fixed simple spectrum $\lambda_1 < \cdots < \lambda_n$
is covered by open dense subsets $\UpLa$, indexed by permutations
$\pi$ and bidiagonal coordinates provide a diffeomorphism between
each $\UpLa$ and $\RR^{n-1}$.
The new
coordinates are used to convert asymptotic issues of Wilkinson's
iteration into local theory in an appropriate chart.

There are two subsets of $\ILa$ in which
Wilkinson's step might break down.
First, let $\Zz \subset \ILa$ be the set of matrices $T$
for which the shift $\omega(T)$ equals an eigenvalue of $T$:
we are especially interested in $\Dd_{0,i} \subset \Zz$,
the set of matrices $T$ with $(T)_{n,n} = \lambda_i$, $(T)_{n,n-1} = 0$.
Second, the two eigenvalues $\omega_+(T) \ge \omega_-(T)$
of the bottom $2 \times 2$ block $\hat T$
may be equally distant from the corner entry $(T)_{n,n}$:
let $\Yy \subset \ILa$ be the set of such matrices $T$.
It is clear that the function $\ww$ is smoothly defined
at least in $\ILa - \Zz - \Yy$ but using bidiagonal coordinates
we shall see that the domain can be taken to be much larger.
Indeed, for a matrix $T \in \UpLa \subset \ILa$
with bidiagonal coordinates $\lambda_i$ and $\beta^\pi_i$,
$T \not\in \Zz \cup \Yy$,
the $\beta^\pi_i$'s of $\ww(T)$ are given by
\[ \ww_\pi(\beta^\pi_1, \ldots, \beta^\pi_{n-1}) =
\left( \left |\frac{\lambda_{\pi(2)} - \omega(T)}{\lambda_{\pi(1)} - \omega(T)} \right | \beta^\pi_1, \ldots,
 \left | \frac{\lambda_{\pi(n)} - \omega(T)}{\lambda_{\pi(n-1)} - \omega(T)} \right | \beta^\pi_{n-1}
\right ). \]
In bidiagonal coordinates, points in $\Dd_{0,\pi(n)} \cap \UpLa$
satisfy $\beta^\pi_{n-1} = 0$ and an inspection of the formula above
shows that $\ww_\pi$ extends smoothly to $\Dd_{0,\pi(n)} - \Yy$.
More generally, near $p_0 \in \Dd_{0,\pi(n)}$
the quotient $\beta^\pi_{n-1}/(T)_{n,n-1}$
is bounded above and below
and rates of convergence are the same in both variables.
We then expand $(\ww_\pi(\beta^\pi_1,\ldots,\beta^\pi_{n-1}))_{n-1}$
in a Taylor series around $p_0$.
Notice that $\omega(T) \approx \lambda_{\pi(n)}$ near $p_0$ and
oddness of this function in the variable $\beta^\pi_{n-1}$ allows for
\[ (\ww_{\pi}(\beta^{\pi}_1, \ldots , \beta^{\pi}_{n-1}))_{n-1} =
G(\beta^{\pi}_1, \ldots , \beta^{\pi}_{n-1}) (\beta^{\pi}_{n-1})^3 \]
for some smooth function $G$, yielding the cubic estimate
\[ |(\ww(T))_{n,n-1}| < C |(T)_{n,n-1}|^3 \]
for some $C > 0$, $T$ in a neighborhood of $p_0$.

Points of $\Yy$ with $\beta^\pi_{n-1} \ne 0$ are step-like discontinuities
for $\ww$.
The behavior of $\ww$ near matrices $p_0 \in \Yy \cap \Zz$
is more complicated;  in figure \ref{fig:142} we show what happens
for $\Lambda = (1,2,4)$, a typical spectrum.
The upshot from the figure is that typically the few points
in $\Yy \cap \Dd_{0,i}$, for which the cubic estimate
does not hold, are isolated in the sequence of iterations $T_k = \ww^k(T_0)$
and are irrelevant in the long run.
More precisely, the argument holds for \textit{AP-free} spectra,
i.e., spectra which contain no three terms arithmetical progression.
In such case, there exist positive constants $C$ and $K$ such that
$|(T_{k+1})_{n,n-1}| > C |(T_k)_{n,n-1}|^3$ holds
for at most $K$ values of $k$ (see theorem \ref{theo:wilk}),
yielding genuine cubic convergence of Wilkinson's iteration.

In the $3 \times 3$ AP case, instead,
there exists a point $p_0 \in \Yy \cap \Dd_{0,i}$
which is kept fixed by $\ww$.
The graphs of $\omega_+$ and $\omega_-$ near $p_0$ resemble
the two branches of the cone $z^2 = xz + y^2$ near the origin,
$x$ and $y$ corresponding to $\beta^\pi_1$ and $\beta^\pi_2$, respectively.
Sections of the cone by planes $x = a$, $a \ne 0$, are hyperbola
with a branch passing through $(a,0,0)$ and $z$ is therefore of the order
of $y^2$ for small $y$: this quadratic behavior of $\omega$ implies
the cubic estimates for the $(n,n-1)$ entry under $\ww$.
For $a = 0$, however, the intersection of the cone with the plane $x = a$
is $z = \pm |y|$; $z$ and therefore $\omega$ are of the order of $|y|$
and $\beta^\pi_2$, respectively:
this entails a {\it quadratic} estimate for the $(n,n-1)$ entry under $\ww$.
Hence, the rate of convergence is dictated by whether $T_k$
remains near $p_0$ when $k$ goes to infinity.
It turns out that $T_k$ tends to $p_0$ (and cubic convergence fails)
if and only if $T_0 \in \Xx$,
where $\Xx$ is a remarkable set (see figure \ref{fig:xxx}).
A point $T_0 \in \Xx$ has a {\it sign sequence} $s: \NN \to \{1,-1\}$,
where $s(k)$ indicates whether $\omega(T_k)$ equals $\omega_+$ or $\omega_-$.
The set $\Ss$ of sign sequences is a Cantor set
and $\Xx$ is the disjoint union of graphs of Lipschitz
functions $f_s: [-a^\ast, a^\ast] \to \RR$
with $f_s(0) = 0$ for all $s \in \Ss$;
a good reference for dynamically defined Cantor sets is \cite{PT}.
For $y_0 \in [-a^\ast, a^\ast]$,
$(f_s(y_0), y_0) \in \Xx$ has sign sequence $s$.
There is an open set of points $T_0 \not\in \Xx$
which are between branches of $\Xx$:
for sufficiently large $k$,
their iterates $T_k$ will be either to the left or to the right of $\Xx$
and cubic convergence applies.
Another application of dynamical systems to numerical spectral theory
is the work of Batterson and Smillie (\cite{BS})
on the Rayleigh quotient iteration.

We begin the paper collecting from \cite{LST} the required
information about bidiagonal coordinates, essentially, the
description of Wilkinson's iteration in terms of bidiagonal
variables. Section 2 also contains a key ingredient: the
construction of a function $h(T)$ which grows along Wilkinson's
iterations. The proof that $h(T)$ indeed satisfies
this property is indirect: it requires the interpretation of a $QR$ step
as the time one map of a Toda flow.
The monotonicity of $h$ then follows by a simple
differentiation argument.

In section 3 we prove cubic convergence of Wilkinson's shift iteration
for AP-free matrices. In sections 4 and 5,
we show that for $3 \times 3$ AP-matrices, the rate of convergence
is usually cubic but strictly quadratic for a thin set of initial
conditions.

The authors acknowledge support from CNPq, CAPES, IM-AGIMB and Faperj.

\section{Preliminaries}

Let $\ILa$ be the set of real, symmetric, tridiagonal matrices with
simple spectrum $\lambda_1,\ldots,\lambda_n$ and set $\Lambda =
\diag(\lambda_1,\ldots,\lambda_n)$. For $T \in \ILa$, write $T =
Q^\transpose \Lambda Q$ for some $Q \in O(n)$. Write the $PLU$
factorization of $Q$, i.e., $Q = PLU$ where $P$ is a permutation
matrix, $L$ is lower unipotent and $U$ is upper triangular. This is
usually possible for $P = P_\pi$ for several permutations $\pi \in
S_n$. Indeed, since $Q$ is invertible, there is a matrix $P^{-1} Q$
obtained by permuting the rows of $Q$ for which all the leading
principal minors are invertible. For a permutation $\pi$, define
$\UpLa$ to be the set of matrices $T \in \ILa$ for which $P_\pi^{-1}
Q$ admits an $LU$ factorization. The following lemma provides other
descriptions of $\UpLa$. Let $\Ee$ be the set of  diagonal matrices
having the values $1$ or $-1$ along the diagonal.

\begin{lemma}[\cite{LST}, lemma 3.1]
\label{lemma:upla}
Take $T \in \ILa$ with unreduced blocks $T_1,\ldots,T_k$ of sizes
$n_1,\ldots,n_k$ along the diagonal.
Then $T \in \UpLa$ if and only if the eigenvalues of $T_i$ are
\[ \lambda^{\pi}_{n_1+\cdots+n_{i-1}+1},\ldots,
\lambda^{\pi}_{n_1+\cdots+n_{i-1}+n_i}. \]
Alternatively,
$\UpLa$ is the union of all open faces in $\ILa$ adjacent to $\Lambda^\pi$.
Also, for $E \in \Ee$, if $T \in \UpLa$ then $ETE \in \UpLa$.
\end{lemma}


For $T \in \UpLa$, we define the {\it $\pi$-normalized
diagonalization} as the unique factorization $T = Q_\pi^\transpose
\Lambda^\pi Q_\pi$ for which the $LU$ factorization of $Q_\pi$
yields a matrix $U$ with positive diagonal. Notice that $Q_{\pi} = E
P_\pi^{-1} Q$ for some $E \in \Ee$.

We now construct charts for $\ILa$, $\phi_\pi: \RR^{n-1} \to \UpLa$ and
$\phi^{-1}_\pi: \UpLa \to \RR^{n-1}$
with $\phi_\pi(0) = \Lambda^\pi \in \UpLa$.
For $T \in \UpLa$, consider its $\pi$-normalized diagonalization
$T = Q_\pi^\transpose \Lambda^\pi Q_\pi$ and $Q = P_\pi Q_\pi$
so that $T = Q^\transpose \Lambda Q$.
Write $Q_\pi = L_\pi U_\pi$ (which is possible because
the leading principal minors of $Q_\pi$ are positive)
and therefore $P_\pi L_\pi = QR_\pi$,
where $R_\pi = U_\pi^{-1}$ is also upper triangular.
Set \[ B_\pi = R_\pi^{-1} T R_\pi = L_\pi^{-1} \Lambda^\pi L_\pi. \]
Notice that the columns of $L_\pi^{-1}$ are the eigenvectors of $B_\pi$.
From $B_\pi = R_\pi^{-1} T R_\pi$, $B_\pi$ is upper Hessenberg and
from $B_\pi = L_\pi^{-1} \Lambda^\pi L_\pi$, it is lower triangular
with diagonal entries $\lambda^\pi_1, \ldots, \lambda^\pi_n$.
Summing up, $B_\pi$ is lower bidiagonal:
\[ B_\pi = \begin{pmatrix}
\lambda^\pi_1 & & & & \\
\beta^\pi_1 & \lambda^\pi_2 & & & \\
& \beta^\pi_2 & \lambda^\pi_3 & & \\
& & \ddots & \ddots & \\
& & & \beta^\pi_{n-1} & \lambda^\pi_n \end{pmatrix}. \] Define
$\psi_\pi$ to be the map just constructed taking $T \in \UpLa$ to
$(\beta^\pi_1,\ldots,\beta^\pi_{n-1})$. We call
$\beta^\pi_1,\ldots,\beta^\pi_{n-1}$ (together with $\lambda^\pi_1,
\ldots, \lambda^\pi_n$) the {\it $\pi$-bidiagonal coordinates} of
$T$.

\begin{theo}[\cite{LST}, theorem 3.4]
\label{theo:bidiagonal} The map $\psi_\pi: \UpLa \to \RR^{n-1}$ is a
diffeomorphism with inverse $\phi_\pi: \RR^{n-1} \to \UpLa$.
\end{theo}

The map $\phi_\pi$ takes open quadrants of $\RR^{n-1}$ diffeomorphically
to the connected components of $\ILa$ formed by irreducible matrices
(connected components are indexed by signs of off-diagonal entries).
Also, the hyperplanes $\beta^\pi_i = 0$ in $\RR^{n-1}$ are taken
diffeomorphically to the set of matrices $T$ in $\UpLa$
with $(T)_{i,i+1} = 0$.
Indeed, it is shown in theorem 4.5 of \cite{LST} that the quotient
$\beta^\pi_i/((T)_{i,i+1})$ is a smooth nonzero function in $\UpLa$.
The following lemma follows by compactness.

\begin{lemma}
\label{lemma:bbeta}
Given a compact subset $K_\pi \subset \UpLa$, there exist positive
constants $C < C'$ such that, for $T \in K_\pi$,
\[ C|(T)_{n,n-1}| \le |\beta^\pi_{n-1}| \le C'|(T)_{n,n-1}|. \]
\end{lemma}

We say that a function $\alpha: \ILa \to \RR$ is {\it even}
if $\alpha(ETE) = \alpha(T)$ for any $T \in \ILa$ and $E \in \Ee$.
The following lemma translates this definition to bidiagonal coordinates.

\begin{lemma}[\cite{LST}, lemma 3.6]
\label{lemma:ETE}
If $E = \diag(\sigma_1, \ldots, \sigma_n) \in \Ee$, $T \in \UpLa$, and
$\psi_\pi(T) = (\beta^\pi_1,\ldots,\beta^\pi_{n-1})$
then
\[ \psi_\pi(ETE) =
(\sigma_1 \sigma_2 \beta^\pi_1,\ldots,
\sigma_{n-1} \sigma_n \beta^\pi_{n-1}). \]

A function $\alpha: \ILa \to \RR$ is even if and only if
each $\alpha \circ \phi_\pi: \RR^{n-1} \to \RR$ is even in each coordinate.
\end{lemma}


As an example of bidiagonal coordinates, let $\Lambda = \diag(-1,0,1)$.
Set $\pi(1) = 3$, $\pi(2) = 1$, $\pi(3) = 2$.
Matrices will be described by their $\pi$-bidiagonal coordinates
$x = \beta_1^\pi$ and $y = \beta_2^\pi$.
Since $B_\pi = L_\pi^{-1} \Lambda^\pi L_\pi$, we obtain
\[
\Lambda^\pi = \begin{pmatrix} 1 & 0 & 0 \\ 0 & -1 & 0 \\
0 & 0 & 0 \end{pmatrix}, \quad
B_\pi = \begin{pmatrix} 1 & 0 & 0 \\ x & -1 & 0 \\
0 & y & 0 \end{pmatrix}, \quad
L_\pi = \begin{pmatrix} 1 & 0 & 0 \\ -x/2 & 1 & 0 \\
- xy & y & 1 \end{pmatrix} \]
and writing $Q_\pi = L_\pi U_\pi$ we have
\[ Q_\pi =
\frac{1}{r_1r_2}
\begin{pmatrix}
{2}{r_2} & {2x(1+2y^2)} & {xy}{r_1} \\
{-x}{r_2} & {2(2+x^2y^2)} & {-2y}{r_1} \\
{-2xy}{r_2} & {y(4-x^2)} & {2}{r_1}
\end{pmatrix}, \]
where $r_1 = \sqrt{4 + x^2 + 4x^2y^2}$ and $r_2 = \sqrt{4 + 4y^2 + x^2y^2}$.
From $T = Q_\pi^\transpose \Lambda^\pi Q_\pi$,
\[ T = \frac{1}{r_1^2r_2^2}
\begin{pmatrix}
(4 - x^2) r_2^2 & 2x r_2^3 & 0 \\
2x r_2^3 & -4(4 - x^2 - 4x^2y^4 + x^4y^4) & 2y r_1^3 \\
0 & 2y r_1^3 & y^2(x^2 - 4) r_1^2
\end{pmatrix}. \]

This example will be revisited in sections 4 and 5.


For an invertible real matrix $M$,
write the unique {\it $QR$ factorization} \[ M = \bQ(M) \bR(M), \]
where $\bQ(M)$ is orthogonal and $\bR(M)$ is upper
triangular with positive diagonal. For some function $f$ taking
nonzero values on the spectrum of a tridiagonal symmetric matrix
$T_0$, the \textit{$QR$ step induced by $f$} is the map
\[ F(T) = \bQ(f(T))^\transpose\,T\,\bQ(f(T)) = \bR(f(T))\,T\,\bR(f(T))^{-1}. \]
From both equalities, we learn that
$F(T)$ is also tridiagonal, with same signs and zeroes along the
off-diagonal entries than $T$. The standard $QR$ step corresponds
to $f(x) = x$ and taking a shift $s$ means taking $f(x) = x-s$.
The map $F$ admits a simple description in terms of bidiagonal coordinates.

\begin{prop}[\cite{LST}, proposition 4.2]
\label{prop:stepbi}
For $f$ taking nonzero values on the spectrum of $T$,
\[ (F \circ \phi_\pi)(\beta^\pi_1, \ldots, \beta^\pi_{n-1}) =
\phi_\pi\left( \left | \frac{f(\lambda_{\pi(2)})}{f(\lambda_{\pi(1)})} \right |
\beta^\pi_1, \ldots,
\left | \frac{f(\lambda_{\pi(n)})}{f(\lambda_{\pi(n-1)})} \right |
\beta^\pi_{n-1} \right ) \]
\end{prop}

We now present some technical results concerning
the dynamics $QR$ type iterations which will be needed in the proof
of theorem \ref{theo:wilk}. 
The argument is rather indirect and seems to require
the language of Toda flows.
For a square matrix $M$, let $S = \Pi_a M$
be the skew symmetric matrix for which $(M)_{ij} = (S)_{ij}$ for $i > j$.
The following result, which follows by direct computation,
relates Toda flows and $QR$ iterations.

\begin{prop}[\cite{Symes2}]
\label{prop:flowiteration}
Let $\Lambda$ be a diagonal matrix with simple spectrum,
$g: \RR \to \RR$ be a smooth function,
$X_g$ be the vector field $X_g(T) = [T, \Pi_a g(T)]$
and $\bT: \RR \to \ILa$ be a path satisfying $\frac{d}{dt}\bT = X_g(\bT)$,
$\bT(0) = T_0$.
Then
\[ \bT(t) = \bQ(\exp(t\,g(T_0)))^\transpose\,T_0\,\bQ(\exp(t\,g(T_0))) =
\bR(\exp(t\,g(T_0))) T_0 \bR(\exp(t\,g(T_0)))^{-1},\]
or, in other words, $\bT(t) = F(T_0)$ where $f(x) = \exp(t\,g(x))$.
\end{prop}

The dynamics of Toda vector fields is rather simple:
they are essentially gradients of Morse functions.
The following result is known for $g(x) = x$
(\cite{BBR}, \cite{DRTW}, \cite{Tomei}).

\begin{theo}
\label{theo:todamorse}
Let $g: \RR \to \RR$ satisfy $g(\lambda_i) \ne g(\lambda_j)$
for distinct eigenvalues $\lambda_i, \lambda_j$ of $\Lambda$.
Let $X_g$ be the Toda vector field $[T, \Pi_a g(T)]$ on $\ILa$.
Let $M = \diag(\mu_1, \mu_2, \ldots, \mu_n)$, $\mu_1 > \mu_2 > \cdots > \mu_n$;
let $h_{M,g}: \ILa \to \RR$ be the smooth function
$h_{M,g}(T) = \trace(M g(T))$.
Then the directional derivative $X_g h_{M,g}$ is strictly positive
except at diagonal matrices.
\end{theo}

\proof
Define a path $\bT: \RR \to \ILa$ satisfying $\frac{d}{dt}\bT = X_g$ as above.
We first claim that $\frac{d}{dt}\tilde g(\bT) = [\tilde g(\bT), \Pi_a g(\bT)]$
for any smooth function $\tilde g$.
Indeed, since Toda flows preserve spectra,
$\tilde g$ may be replaced by a polynomial $p$
which coincides with $\tilde g$ on the spectrum of $\Lambda$.
By linearity, it suffices to consider $p_k(x) = x^k$,
which is handled by induction on $k$:
\begin{align}
\frac{d}{dt}p_{k+1}(\bT) &= \frac{d}{dt} (\bT p_k(\bT))
= \bT \left(\frac{d}{dt} p_k(\bT) \right) +
\left(\frac{d}{dt} \bT \right) p_k(\bT) \notag \\
&= \bT [\bT^k, \Pi_a f(\bT)] + [\bT, \Pi_a f(\bT)] \bT^k \notag \\
&= \bT^{k+1} (\Pi_a f(\bT)) - \bT (\Pi_a f(\bT)) \bT^k +
\bT (\Pi_a f(\bT)) \bT^k - (\Pi_a f(\bT)) \bT^{k+1} \notag \\
&= [p_{k+1}(\bT), \Pi_a f(\bT)]. \notag
\end{align}
Take $\tilde g = g$ and compute
the derivative of $h_{M,g}$ along the path $\bT$:
\begin{align}
X_g h_{M,g} &= \frac{d}{dt} h_{M,g}(\bT) =
\frac{d}{dt} \left( \trace\; M g(\bT) \right) \notag\\
&= \trace\; M \left(\frac{d}{dt} g(\bT) \right)
= \trace (M [g(\bT), \Pi_a g(\bT)]) \notag\\
&= \sum_{1 \le i < j \le n} 2(\mu_i - \mu_j) (g(\bT))_{ij}^2. \notag
\end{align}
Since $g(T)$ has simple spectrum,
it is diagonal only if $T$ also is and we are done.
\qed

\begin{coro}
\label{coro:QRmorse}
Let $f: \RR \to \RR$ be a function satisfying
$|f(\lambda_i)| \ne |f(\lambda_j)| \ne 0$ for $i \ne j$
and let $F: \ILa \to \ILa$ be the $QR$ step induced by $f$.
Let $\Kk \subset \ILa$ be a compact set containing no diagonal matrices.
Then there exists $K > 0$ such that, for all $T \in \ILa$,
\[ |\{k \in \NN | F^k(T) \in \Kk \}| < K. \]
\end{coro}

\proof
Consider $g(x) = \log|f(x)|$, $X_g$, $\bT$, $M$ and $h_{M,g}$
as in the previous theorem.
Since $\Kk$ is compact and avoids diagonal matrices,
$X_g h_{M,g} > \epsilon > 0$ on $\Kk$.
From proposition \ref{prop:flowiteration},
$F^k(T_0) = \bT(k)$ and we are done.
\qed

It is not true that, given $\Kk$,
there exists $K_2$ such that, for all $T_0 \in \ILa$,
$F^k(T_0) \not\in \Kk$ for all $k > K_2$.
For instance, the exit time from a small neighborhood
of a diagonal matrix $T$ is not uniformly bounded;
another situation is given in figure \ref{fig:notsocubic} below.

\section{Wilkinson's shift and AP-free matrices}

We consider iterations of shifted $QR$ steps, induced by $f_s(x) =
x-s$. The function $F_s: \ILa \to \ILa$ is defined whenever $s$ is
not an eigenvalue. Frequently, the shift $s$ is taken to depend on
$T$. We consider the asymptotic properties of a special choice of
shift, which defines {\it Wilkinson's step} (\cite{Wilkinson}). For
$T \in \ILa$, let $\omega_+(T) \ge \omega_-(T)$ be the two
eigenvalues of the bottom $2 \times 2$ block of $T$ and set
$\omega(T)$ to be the one nearest to the entry $(T)_{n,n}$. This
defines continuous maps $\omega_+, \omega_-: \ILa \to \RR$ which are
smooth in $\ILa - \Yy_0$, where $\Yy_0 \subset \ILa$ is the set of
matrices $T$ for which $(T)_{n,n} = (T)_{n-1,n-1}$ and $(T)_{n,n-1}
= 0$. Indeed, $T \in \Yy_0$ if and only if $\omega_+(T) =
\omega_-(T)$. Also, $\omega: \ILa - \Yy \to \RR$ is smooth where $T
\in \Yy \subset \ILa$ if $(T)_{n,n} = (T)_{n-1,n-1}$. Indeed, $T \in
\Yy$ if and only if $(T)_{n,n}$ is equally distant from both
eigenvalues of its bottom $2 \times 2$ block, which is the only way
smoothness could fail. Also, $\omega_\pm$ and $\omega$ are even in
the sense of lemma \ref{lemma:ETE}.

\begin{lemma}
\label{lemma:wlip}
The functions $\omega_\pm: \ILa \to \RR$ are Lipschitz.
\end{lemma}

\proof
The function taking $T$ to its bottom $2\times 2$ block is clearly smooth.
The function taking a $2\times 2$ symmetric matrix $A$ to its larger
(resp. smaller) eigenvalue is Lipschitz in compact sets.
The lemma follows from composition.
\qed

Let $ \Zz_k = \{ T \in \ILa;\; \omega(T) = \lambda_k \} $,
$\Zz = \bigcup_k \Zz_k$, $\hat\Zz_k = \Zz - \Zz_k$.
Notice that if $(T)_{n,n-1} = 0$ then $\omega(T) = (T)_{n,n} = \lambda_k$
for some $k$ and therefore $T \in \Zz_k$.
Also, if $(T)_{n-1,n-2} = 0$ then again $\omega(T) = \lambda_k$
for some $k$.
Define $\ww: \ILa - \Yy - \Zz \to \ILa$ by $\ww(T) = F_{\omega(T)}(T)$
for $f_s(x) = x-s$. Notice that $\ww$ is an odd function.

\begin{figure}[ht]
\begin{center}
\psfrag{Yy}{$\Yy$}
\epsfig{height=28mm,file=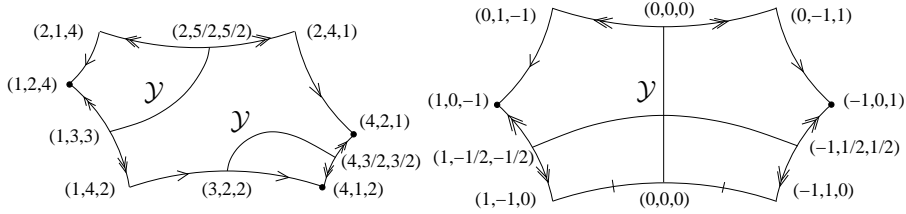}
\end{center}
\caption{\capsize The phase space of Wilkinson's step for $n=3$.}
\label{fig:wil}
\end{figure}

Figure \ref{fig:wil} shows $\JLa$ for $\Lambda = \diag(1,2,4)$
and $\Lambda = \diag(-1,0,1)$.
Labels indicate the diagonal entries and vertices are diagonal matrices.
The set $\Yy$, on which $\ww$ is not defined,
degenerates in the second example.
Vertices are fixed points and edges are invariant under $\ww$.
A simple arrow indicates the order of the points $T, \ww(T), \ww^2(T), \ldots$
along the edge.
For points $T$ on an arc labeled with a double arrow, $\ww(T)$ is a vertex.
Arcs marked with a transversal segment consist of fixed points of $\ww$.

We apply proposition \ref{prop:stepbi} to write down Wilkinson's step
in bidiagonal coordinates.
Define $\ww_\pi$ by
\[ \ww_\pi(\beta^\pi_1, \ldots, \beta^\pi_{n-1}) =
\left( \left |
\frac{\lambda_{\pi(2)} - \omega}{\lambda_{\pi(1)} - \omega}
\right |  \beta^\pi_1, \ldots, \left |
\frac{\lambda_{\pi(n)} - \omega}{\lambda_{\pi(n-1)} - \omega}
\right | \beta^\pi_{n-1} \right ) \]
where $\omega = \omega(\phi_\pi(\beta^\pi_1, \ldots, \beta^\pi_{n-1}))$.
Thus, the natural domain for $\ww_\pi$ is
$\RR^{n-1} - \phi_\pi^{-1}(\Yy \cap \UpLa) -
\phi_\pi^{-1}(\hat\Zz_{\pi(n)} \cap \UpLa)$,
where $\ww_\pi$ is a smooth function,
indicating that points in $\phi_\pi^{-1}(\Zz_{\pi(n)})$ are
removable singularities and that, despite absolute values in the formula,
$\ww_\pi$ is smooth at such points.
Notice that $\ww_\pi$ is odd, since $\omega$ is even
in each variable $\beta^\pi_i$.
Also, points in $\phi_\pi^{-1}(\Zz_{\pi(n)})$ are of the form
$\beta^\pi_{n-1} = 0$ or, equivalently, $(T)_{n,n-1} = 0$.

\begin{figure}[ht]
\begin{center}
\psfrag{D0}{}
\psfrag{Ddm}{$\Dd_-$}
\psfrag{Ddp}{$\Dd_+$}
\psfrag{Yy}{$\Yy$}
\psfrag{T0}{$p_0$}
\psfrag{wT0}{$\ww(p_0)$}
\epsfig{height=50mm,file=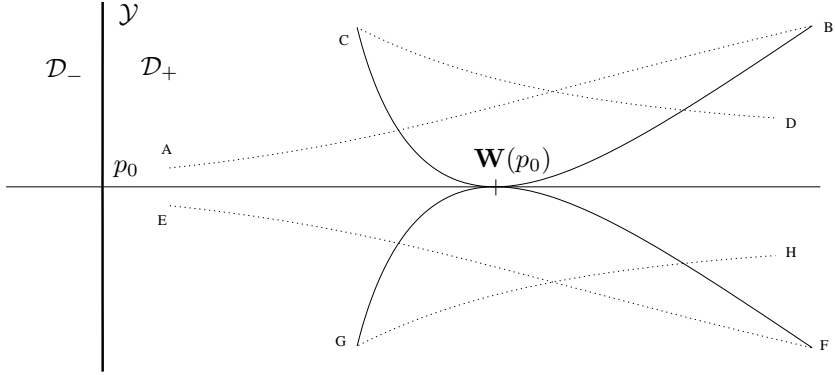}
\end{center}
\caption{\capsize Image of $\Yy$ under $\ww_\pi$ for $\Lambda = \diag(1,2,4)$
(stretched vertically).}
\label{fig:142}
\end{figure}

Points of $\Yy$ with $\beta^\pi_{n-1} \ne 0$ are step-like discontinuities
for $\ww$.
The behavior of $\ww$ near matrices $p_0 \in \Yy \cap \Zz$
is more complicated;  in figure \ref{fig:142} we show what happens
for $\Lambda = (1,2,4)$, a typical $3 \times 3$ AP-free spectrum.
Close to
\[ p_0 = \begin{pmatrix} 3 & \sqrt{2} & 0 \\ \sqrt{2} & 2 & 0 \\
0 & 0 & 2 \end{pmatrix}, \quad \pi(1) = 1, \pi(2) = 3, \pi(3) = 2, \quad
\beta^\pi_1 = 3 \sqrt{2}, \beta^\pi_2 = 0,\]
the set $\Yy$ divides the plane into two sides $\Dd_+$ and $\Dd_-$,
where we take $\omega_+$ and $\omega_-$, respectively,
in the definition of $\ww$.
From each side $\Dd_\pm$, the function $\ww$ can be continuously extended
to $\Yy$ but the two values thus obtained are quite different
except at $p_0$.

The figure is drawn with a vertical stretching factor of $200$.
The thick vertical curve (which looks like a straight line due to stretching)
is the set $\Yy$.
To the right, the curve $BF$ (with a cusp at $\ww(p_0)$) is the image of
the arc in $\Yy$ from $\beta^\pi_2 = 0.1$ to $\beta^\pi_2 = -0.1$
under $\ww$ with the choice $\omega = \omega_-$;
to the left, the curve $CG$ is the image of the same arc,
now with $\omega = \omega_+$.
The dotted lines $AB$ and $CD$ are the image under $\ww$ of the horizontal line
$\beta^\pi_2 = 0.1$: notice the jump discontinuity at $\Yy$ from $B$ to $C$.
Oddness of $\ww$ in the coordinate $\beta^\pi_2$ explains the mirror symmetry
in the horizontal axis.

The study of the asymptotic behavior of the iterations of
Wilkinson's step becomes significantly simpler
by taking into account the following result.

\begin{theo}[\cite{Parlett}, p. 152]
\label{theo:parlett}
If $T_k = \ww^k(T_0)$ then $\lim_{k \to \infty}(T_k)_{n,n-1} = 0$.
\end{theo}

From this result, deflation is always possible in numerical
implementations of Wilkinson's iteration: $T$ will be truncated
so as to split into blocks of size $n-1$ and $1$.
We will investigate the rate of convergence of $(T_k)_{n,n-1}$ to $0$.

\begin{prop}
\label{prop:DqDc}
Let $\Dd(\epsilon) = \{ T \in \ILa \; | \; |(T)_{n,n-1}| < \epsilon \}$.
Given a spectrum $\lambda_1 < \cdots < \lambda_n$,
there exists $\epsilon_q > 0$ and $C_q > 0$ such that
$\ww(\Dd(\epsilon)) \subset \Dd(\epsilon)$ for any $\epsilon \le \epsilon_q$
and $ |(\ww(T))_{n,n-1}| \le C_q |(T)_{n,n-1}|^2$ for $T \in \Dd(\epsilon_q)$.

Furthermore, for any subset $\Dd_c \subset \Dd(\epsilon_q)$ satisfying
$\overline{\Dd_c} \cap \Yy_0 = \emptyset$
there exists $C_c > 0$ such that
$|(\ww(T))_{n,n-1}| \le C_c |(T)_{n,n-1}|^3$ for $T \in \Dd_c$.
\end{prop}

\proof
Define smaller open sets $\Vv^\pi \subset K_\pi \subset \UpLa$
where each $K_\pi$ is compact such that the open sets $\Vv^\pi$
still cover $\ILa$.
Let $\Dd_0 \subset \ILa$ be the set of matrices $T$ with $(T)_{n,n-1} = 0$.
The set $\Dd_0$ is a compact submanifold of codimension $1$.
The set $\Dd_0$ is not contained in any $\Vv^\pi$ but is clearly
covered by them and we can then write the $\pi$-bidiagonal coordinates
$\beta^\pi_1, \cdots, \beta^\pi_{n-1}$ for $T \in \Vv^\pi$.
Lemma \ref{lemma:bbeta} allows us to identify in $\Vv^\pi$
rates of decay for $(T)_{n,n-1}$ and $\beta^\pi_{n-1}$.

Let $\epsilon_1 < |\lambda_i - \lambda_j|/2$ for all $i \ne j$.
We can take $\epsilon_0 > 0$ such that $T \in \Dd(\epsilon_0)$
implies $|(T)_{n,n} - \lambda_i| < \epsilon_1$ for some (unique) $i$
($i$ depends on $T$).
Call $\Dd^i(\epsilon_0) \subset \Dd(\epsilon_0)$ the set of such $T$.
Assume without loss that $\Dd^i(\epsilon_0)$ is covered by $\Vv^\pi$ for
permutations $\pi$ satisfying $\pi(n) = i$.
In $\Dd^i(\epsilon_0)$, $\ww$ in $\pi$-bidiagonal coordinates
is given by
\[ \ww_\pi(\beta^\pi_1, \ldots, \beta^\pi_{n-1}) =
\left( \left | \frac{\lambda_{\pi(2)} - \omega}{\lambda_{\pi(1)} - \omega} \right | \beta^\pi_1,
\ldots,
\left | \frac{\lambda_{i} - \omega}{\lambda_{\pi(n-1)} - \omega} \right | \beta^\pi_{n-1}
\right ). \]
In each $\Vv^\pi \cap \Dd^i(\epsilon_0)$, $\omega$ is a continuous
function taking the value $\lambda_i$ when $\beta^\pi_{n-1} = 0$.
Since $\omega_\pm$ are Lipschitz functions of $T$,
$|\lambda_i - \omega| < L |\beta^\pi_{n-1}|$ for some $L > 0$,
implying from the formula for $\ww_\pi$
the quadratic rate of decay for $\beta^\pi_{n-1}$.
By taking $\epsilon_0$ even smaller,
we can assume that $|\lambda_i - \omega| < \epsilon_1$ in this set
so that the quotients
\[ \frac{\lambda_{\pi(2)} - \omega}{\lambda_{\pi(1)} - \omega}, \cdots,
\frac{\lambda_{\pi(n-1)} - \omega}{\lambda_{\pi(n-2)} - \omega} \]
have absolute value bounded and bounded away from zero
and the first claim follows by compactness.

The Taylor expansion of the last coordinate
of $\ww_{\pi}$ centered at points $T \in (\Vv^\pi \cap \Dd_0) - \Yy_0$
with respect to the variable $\beta^{\pi}_{n-1}$ is of the form
\[ (\ww_{\pi}(\beta^{\pi}_1, \ldots , \beta^{\pi}_{n-1}))_{n-1} =
\sum_\ell a_\ell(\beta^{\pi}_1, \ldots , \beta^{\pi}_{n-2})
(\beta^{\pi}_{n-1})^\ell. \]
The function $\ww_{\pi}$ is odd, so $a_0 = a_2 = 0$.
The factor $(\lambda_{\pi(n)} - \omega)/(\lambda_{\pi(n-1)} - \omega)$
equals $0$ if $\beta^{\pi}_{n-1} = 0$ and therefore $a_1 = 0$.
In other words,
\[ (\ww_{\pi}(\beta^{\pi}_1, \ldots , \beta^{\pi}_{n-1}))_{n-1} =
G(\beta^{\pi}_1, \ldots , \beta^{\pi}_{n-1}) (\beta^{\pi}_{n-1})^3 \]
for some real analytic function $G$.
Again by compactness, we are done.
\qed

A matrix $T$ is an {\it AP-matrix} if its spectrum
contains an arithmetic progression with $3$ terms,
i.e., some eigenvalue is the average of two others;
otherwise, $T$ is {\it AP-free}.
We also refer to AP-free spectra, Jacobi cells, isospectral manifolds
and so on, with the obvious meanings.
The left hexagon in figure \ref{fig:wil} is AP-free and
the right hexagon is not.

\begin{theo}
\label{theo:wilk}
For AP-free tridiagonal matrices, Wilkinson's iteration
has cubic convergence. More precisely, given an AP-free
spectrum $\lambda_1 < \cdots < \lambda_n$,
there exist $C > 0$ and $K > 0$ such that, for any $T_0 \in \ILa$,
\begin{enumerate}[\rm (a)]
\item{if $k > K$ then
$|(T_{k+1})_{n,n-1}| < 1/C$ and
$|(T_{k+1})_{n,n-1}| \le C |(T_k)_{n,n-1}|^2$; }
\item{the set of positive integers $k$ for which
$|(T_{k+1})_{n,n-1}| > C |(T_k)_{n,n-1}|^3$
has at most $K$ elements.}
\end{enumerate}
\end{theo}

\proof
We keep the notation of the proof of proposition \ref{prop:DqDc}.
Item (a) follows from theorem \ref{theo:parlett}, proposition \ref{prop:DqDc}
and compactness.
Indeed,
for any given $\epsilon > 0$ there exists $K_1$ such that
for any $T_0 \in \ILa$ and for any $k > K_1$
we have $T_k \in \Dd(\epsilon)$.
For item (b),
we need to prove that, given an open $\Dd_c \subset \ILa$
containing the diagonal matrices,
there exists $K_2$ such that, for any $T_0 \in \ILa$,
the set of positive integers $k$ for which
$T_k \not\in \Dd_c$ has at most $K_2$ elements.
Notice that $\Yy_0$ is removed from diagonal matrices
and therefore such a set $\Dd_c$ exists.

As a warm-up case,
let $\Dd_{0,i} \subset \Zz_i$ be the set of matrices $T \in \ILa$
for which $(T)_{n,n} = \lambda_i$, $(T)_{n,n-1} = 0$.
Clearly, $\Dd_{0,i}$ are the $n$ connected components of $\Dd_0$ and
the set $\Dd_{0,i}$ is diffeomorphic to $\ILai$, where
\[ \Lambda_i = \diag(\lambda_1, \ldots, \lambda_{i-1},
\lambda_{i+1}, \ldots \lambda_n), \]
that is, $\Lambda$ without $\lambda_i$.
Points in $\Dd_{0,i}$ are removable singularities of $\ww$ and there
$\ww$ is a $QR$ step on $\ILai$ (the top $(n-1) \times (n-1)$ block)
with $f(x) = x - \lambda_i$.
Apply corollary \ref{coro:QRmorse} to conclude this special case:
simplicity of $\Lambda$ ensures that $\lambda_i$ is not an eigenvalue
of the top $(n-1) \times (n-1)$ block and the fact that $\Lambda$ is AP-free
ensures that $|lambda_j - \lambda_i| \ne |\lambda_{j'} - \lambda_i|$ for 
distinct eigenvalues $\lambda_j$ and $\lambda_{j'}$ of said block.
In general, we need to extend the reasoning in corollary \ref{coro:QRmorse}
to a neighborhood of $\Dd_{0,i}$.

Let $b: \ILa \to \RR$ be the smooth function defined by $b(T) = (T)_{n,n-1}$.
From lemma \ref{lemma:bbeta}, $b$ is transversal to $\Dd_0$;
let $Z$ be a smooth vector field in $\ILa$ such that
the directional derivative $Z b$ satisfies $Z b = 1$
in $\Dd(\epsilon_Z)$ for some $\epsilon_Z > 0$, $\epsilon_Z < \epsilon_q$.
The vector field $Z$ can be integrated to yield a diffeomorphism
$\zeta: \Dd_0 \times (-\epsilon_Z, \epsilon_Z) \to \Dd(\epsilon_Z)$.
Let $\Dd_i(\epsilon_Z) = \zeta(\Dd_{0,i} \times (-\epsilon_Z, \epsilon_Z))$:
we have a diffeomorphism $\tilde\zeta: \ILaie \to \Dd_i(\epsilon_Z)$.
Set $\hat\ww_i: \ILaie \to \ILaie$,
$\hat\ww_i = (\tilde\zeta)^{-1} \circ \ww \circ \tilde\zeta$ and,
for $\ttt_0 \in \ILaie$,
$\ttt_k$ is defined recursively by $\ttt_{k+1} = \hat\ww_i(\ttt_k)$.

Let $g(x) = \log(|x-\lambda_i|)$ and $M = \diag(n-1, n-2, \ldots, 2, 1)$.
Let $X_g$ tangent to $\ILai$ and $h_{M,g}: \ILai \to \RR$ be
as in theorem \ref{theo:todamorse}.
Extend $h_{M,g}$ to $h_{M,g}: \ILaie \to \RR$
by ignoring the second coordinate.
Let $\bT: \RR \to \ILai$ with $\frac{d}{dt}\bT = X_g$:
from proposition \ref{prop:flowiteration}
we have $\hat\ww_i(\bT(0),0) = (\bT(1),0)$.
Therefore, from theorem \ref{theo:todamorse},
$h_{M,g}(\hat\ww_i(\ttt_0)) > h_{M,g}(\ttt_0)$
for any $\ttt_0 = (T_0,0) \in \ILai \times \{0\}$,
$T_0$ not a diagonal matrix.
Set $\delta: \ILaie \to \RR$,
$\delta(\ttt_0) = h_{M,g}(\hat\ww_i(\ttt_0)) - h_{M,g}(\ttt_0)$;
by compactness of the complement of $\Dd_c$,
there exists $\epsilon > 0$ such that
$\delta(\ttt_0) > \epsilon$
if $\ttt_0 \in (\ILai \times \{0\}) - (\tilde\zeta^{-1}(\Dd_c))$.
Since $\delta$ is continuous in $\ILai \times \{0\}$,
there exists $\epsilon' > 0$, $\epsilon' < \epsilon_Z$,
such that $\ttt_0 \in (\ILaiei)  - \tilde\zeta^{-1}(\Dd_c)$
implies $\delta(\ttt_0) > \epsilon/2$.
Clearly, for any $\ttt_0$ in $\ILaiei$,
\[ \sum_{k \ge 0} \delta(\ttt_k) \le \max h_{M,g} - \min h_{M,g} \]
and therefore there are at most
$2(\max h_{M,g} - \min h_{M,g})/\epsilon$ values of $k$
for which $\ttt_k \not\in \tilde\zeta^{-1}(\Dd_c)$
and we are done.
\qed

\begin{figure}[ht]
\begin{center}
\psfrag{Yy}{$\Yy$}
\psfrag{S}{$S$}
\epsfig{height=28mm,file=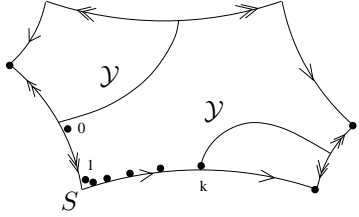}
\end{center}
\caption{\capsize We may have $T_k \in \Yy$ for large values of $k$.}
\label{fig:notsocubic}
\end{figure}

On the other hand, it is {\it not} true that given an AP-free spectrum
$\Lambda$ there exist $C > 0$ and $K$ such that
$|(T_{k+1})_{n,n-1}| \le C |(T_k)_{n,n-1}|^3$ for all $k > K$.
A counterexample is indicated in figure \ref{fig:notsocubic}:
the orbit may spend an arbitrarily large number of steps near
the saddle point $S$ and we may therefore have $T_k \in \Yy$
for arbitrarily large $k$.



\section{Wilkinson's step for $3\times 3$ AP-matrices}

In this and the following sections
we prove that for any spectrum of the form
$\{a-b, a, a+b\}$ there exist matrices $T_0$ for which
\[ \lim_{k \to \infty} T_k =
\begin{pmatrix} a & b & 0 \\ b & a & 0 \\ 0 & 0 & a \end{pmatrix} \]
and that the convergence of $(T_k)_{n,n-1}$ towards $0$
is still quadratic but not cubic.
For $a = 0$, $b = 1$, this limit is the point
labeled $(0,0,0)$ at the bottom of  the right hexagon
in figure \ref{fig:wil}.

Up to normalizations, a $3\times 3$ AP-matrix
is isospectral to $\Lambda = \diag(-1,0,1)$.
We use the $\pi$-bidiagonal coordinates for matrices in $\UpLa$
computed in section 2 for an appropriate permutation $\pi$.
The bottom $2 \times 2$ block of $T \in \UpLa$ is
\[ \hat T = \frac{1}{r_1^2r_2^2}
\begin{pmatrix}
 -4(x^2 - 4)(x^2y^4 - 1) & 2y r_1^3 \\
 2y r_1^3 & y^2(x^2 - 4) r_1^2
\end{pmatrix} \]
where $r_1 = \sqrt{4+x^2+4x^2y^2}$, $r_2 = \sqrt{4+4y^2+x^2y^2}$.
Let $\omega_+ \ge \omega_-$ be the (real) eigenvalues of $\hat T$;
by interlacing, $-1 \le \omega_- \le 0 \le \omega_+ \le 1$,
with equality only when $x = 0$ or $y = 0$.
The discriminant of the characteristic polynomial of $\hat T$ is
$\Delta = ((x+2)^2 + 8x^2y^2)((x-2)^2 + 8x^2y^2) \ge 0$
which is zero exactly at the points $\pm p_0$, $p_0 = (2, 0)$.
These points correspond to the matrices $\pm P_0$,
\[ P_0 = \begin{pmatrix} 0 & 1 & 0 \\ 1 & 0 & 0 \\ 0 & 0 & 0 \end{pmatrix} \]
and the origin to the matrix $\Lambda^\pi$.

\begin{figure}[ht]
\begin{center}
\psfrag{Yy}{$\Yy$}
\psfrag{omegaplus}{$\omega_+$}
\psfrag{omegaminus}{$\omega_-$}
\psfrag{x}{$x$}
\psfrag{y}{$y$}
\psfrag{Rp}{$R_+$}
\psfrag{Rm}{$R_-$}
\psfrag{r}{$r$}
\psfrag{(0,1)}{$(0,1)$}
\psfrag{(2,0)}{$(2,0)$}
\psfrag{(2,1/2)}{$(2,1/2)$}
\epsfig{height=35mm,file=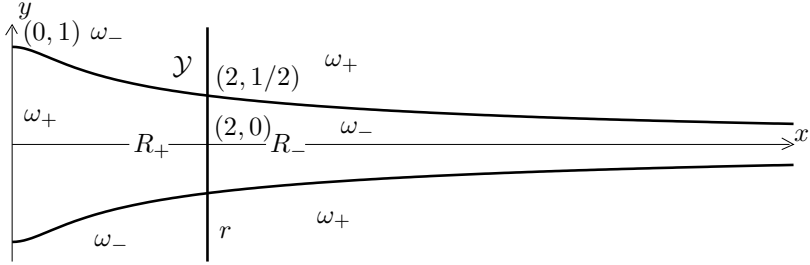}
\end{center}
\caption{\capsize The choice of $\omega$.}
\label{fig:xy}
\end{figure}

The eigenvalue closest to the bottom element $(\hat T)_{2,2} = (T)_{3,3}$
is $\omega_+$ if and only if $(T)_{3,3} > (\omega_+ + \omega_-)/2$.
A straightforward computation yields
\[ (T)_{3,3} - \frac{\omega_+ + \omega_-}{2} =
\frac{(2-x)(2+x)(4 - 4y^2 - x^2y^2 - 8x^2y^4)}{2 r_1^2 r_2^2} \]
which indicates that the $xy$-plane is divided by $\Yy$
into regions, the choice between $\omega_+$ and $\omega_-$
being as in figure \ref{fig:xy}
(where only the region $x > 0$ is shown).
The hexagon of matrices on the right of figure \ref{fig:wil}
is the image under $\phi_\pi$ of the first quadrant
(of $\pi$-bidiagonal coordinates) in figure \ref{fig:xy};
the reader may check, for instance, that $\phi_\pi(2,0) = P_0$
and that
\[ \phi_\pi(0,1) = 
\begin{pmatrix} 1 & 0 & 0 \\ 0 & -1/2 & 1/2 \\ 0 & 1/2 & -1/2 \end{pmatrix}. \]

We consider the region
\[ R = \{ (x,y) | x > 0, 4 - 4y^2 - x^2y^2 - 8x^2y^4 \ge 0 \} \]
and the subsets $R_+ = R \cap ((0,2] \times \RR)$,
$R_- = R \cap ([2,+\infty) \times \RR)$.
For $0 < a \le 1/10$, define the {\it wedge of height $a$} to be
\[ V_a = \{(x,y) | |y| \le a, |y| \ge |x-2|/10 \}. \]

\begin{lemma}
\label{lemma:omega}
The functions $\omega_\pm$ are smooth in $(x,y) \in R$
except at the point $p_0$, where they have a cone-like behavior:
\[ \omega_\pm = \frac{(x-2) \pm \sqrt{(x-2)^2 + 32y^2}}{4} +
O((x-2)^2 + 32y^2). \]
For $(x,y) \in R_+$ (resp. $R_-$),
$0 \le \omega_+ \le 2|y|$ (resp. $-2|y| \le \omega_- \le 0$).
There exists a positive constant $C$ such that
$|\omega_\pm| \ge C |y|$ for $(x,y) \in V_a$.
\end{lemma}

\proof
We have
\[ \omega_\pm = \frac{-4 + x^2 \pm \sqrt\Delta}{2r_1^2}. \]
The displayed estimate for $\omega_\pm$ follows directly from
\[ \lim_{(x,y) \to (2,0)} \frac{\Delta}{16((x - 2)^2 + 32y^2)} = 1.\]
The signs of $\omega_\pm$ follow from interlacing
and the other estimates are now easy.
\qed

\begin{lemma}
\label{lemma:signofomegax}
The partial derivatives $(\omega_\pm)_x$ and $(\omega_\pm)_y$
are uniformly bounded in $R - \{p_0\}$.
For all $(x,y) \in R_\pm - \{p_0\}$ we have $(\omega_\pm)_x \ge 0$,
with equality exactly when $y = 0$.
Also, $(\omega_\pm)_x > 1/120$ in any wedge $V_a$.
Furthermore, for $y \ne 0$, $\pm y (\omega_\pm)_y > 0$.
\end{lemma}

\proof
A straightforward computation yields
\begin{align}
(\omega_\pm)_x &=
\frac{8x}{r_1^4 \sqrt{\Delta}}
\left( \left( (1+2y^2) \sqrt{\Delta} \right) \pm
\left(-4 + x^2 + 8y^2 + 6x^2y^2 + 16x^2y^4 \right) \right), \notag\\
(\omega_\pm)_y &=
\frac{4x^2y}{r_1^4 \sqrt{\Delta}}
\left( \left( (4-x^2) \sqrt{\Delta} \right) \pm
\left(16 + 24x^2 + x^4 + 32x^2y^2 + 8x^4y^2 \right) \right). \notag
\end{align}
Also,
\[ \left( (1+2y^2) \sqrt{\Delta} \right)^2 -
\left(-4 + x^2 + 8y^2 + 6x^2y^2 + 16x^2y^4 \right)^2
= 8 y^2 r_1^4 \ge 0 \]
whence
\[ (1+2y^2) \sqrt{\Delta} \ge
\left|-4 + x^2 + 8y^2 + 6x^2y^2 + 16x^2y^4 \right|, \]
where the equality holds if and only if $y = 0$.
In order to prove the estimate in $V$, write
\begin{align}
(\omega_\pm)_x &=
\frac{8x}{r_1^4 \sqrt{\Delta}}\; \frac{8 y^2 r_1^4}%
{\left( \left( (1+2y^2) \sqrt{\Delta} \right) \mp
\left(-4 + x^2 + 8y^2 + 6x^2y^2 + 16x^2y^4 \right) \right)} \notag\\
&\ge \frac{32 x}{(1+2y^2)((x+2)^2 + 8x^2y^2)}\; \frac{y^2}{(x-2)^2 + 8x^2y^2}
> 1/120. \notag
\end{align}

Since
\[ \left( (4-x^2) \sqrt{\Delta} \right)^2 -
\left(16 + 24x^2 + x^4 + 32x^2y^2 + 8x^4y^2 \right)^2 =
-64x^2 r_1^4 \le 2 \]
we have
\[ \left| (4-x^2) \sqrt{\Delta} \right| \le
16 + 24x^2 + x^4 + 32x^2y^2 + 8x^4y^2, \]
which yields the sign of $(\omega_\pm)_y$.

Boundedness near $x = \infty$ follows from the rates in $x$
since $y$ is bounded.
For $(x,y)$ near $p_0$, expand the formula for $(\omega_\pm)_x$
as a sum of two terms.
The first term is bounded since $\sqrt{\Delta}$ simplifies;
from the computations above, the absolute value of the second term
is no larger.
Since $y/\sqrt{\Delta}$ is bounded near $p_0$, so is $(\omega_\pm)_y$.
\qed

\begin{figure}[ht]
\begin{center}
\epsfig{height=30mm,file=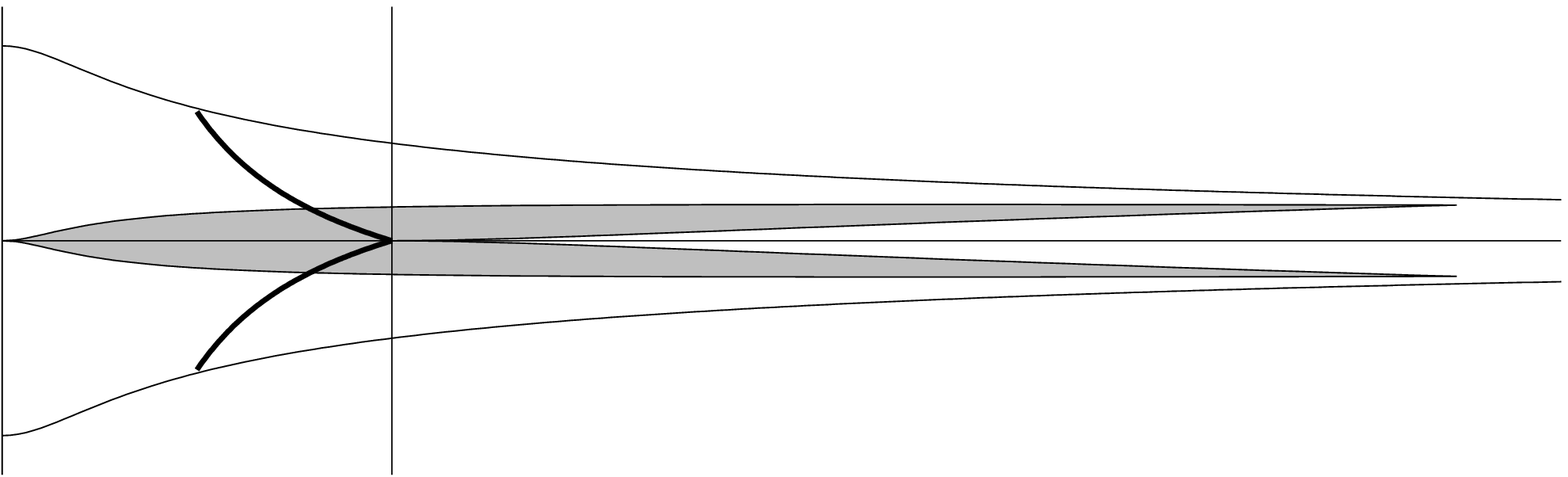}
\end{center}
\begin{center}
\epsfig{height=30mm,file=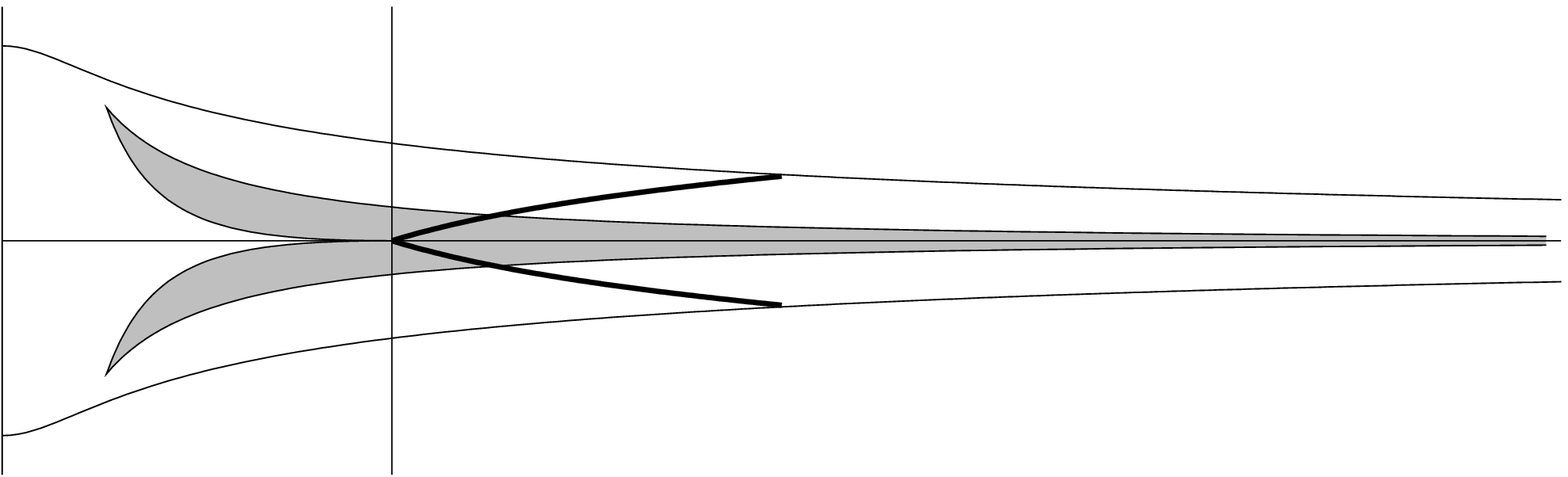}
\end{center}
\caption{\capsize $\ww_\pm(R_\pm)$ (shaded)
and $\ww_\pm^{-1}(r)$ (thick); in scale.}
\label{fig:wrpm}
\end{figure}

In bidiagonal coordinates, Wilkinson's step is given by
\[ \ww(x,y) = \left( \frac{1+\omega}{1-\omega} x,
\frac{|\omega|}{1+\omega} y \right). \]
Since from interlacing $\omega_\pm$ does not change sign, the restrictions
$\ww_+: R_+ \to R \subset \RR^2$ and $\ww_-: R_- \to R \subset \RR^2$,
are continuous in their respective domains and smooth except at $p_0$.
The restrictions of $\ww$ to the the left and to the right
of the vertical line $r$ given by $x = 2$
coincide with $\ww_+$ and $\ww_-$ and the restrictions
of these two functions to $r$ yield different values.
Figure \ref{fig:wrpm} shows
the images of $\ww_+$ and $\ww_-$, clearly contained in $R$.
As we shall see, both $\ww_+$ and $\ww_-$ are homeomorphisms onto
their respective images.

\begin{figure}[ht]
\begin{center}
\psfrag{wmr}{$\ww_-(r)$}
\psfrag{wpr}{$\ww_+(r)$}
\psfrag{wmir}{$\ww_-^{-1}(r)$}
\psfrag{wpir}{$\ww_+^{-1}(r)$}
\psfrag{X}{$p_0$}
\psfrag{r}{$r$}
\epsfig{height=50mm,file=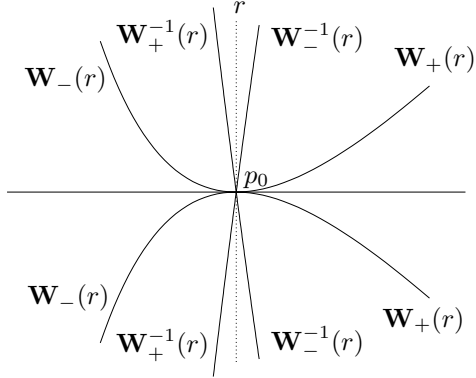}
\end{center}
\caption{\capsize Images and preimages of the line $x=2$; not in scale. }
\label{fig:r5}
\end{figure}

The line $r$ is taken by $\ww_+$ (resp. $\ww_-$) to an arc
contained in $R_-$ (resp. $R_+$), with a cusp at $p_0$.
The horizontal axis is a common tangent to the four smooth subarcs
in the images of $r$.
A straightforward computation verifies that the preimage of the
vertical line $r$ under $\ww$ consists of the two smooth arcs
\[ \left(x, \pm\frac{(x-2)\sqrt{x(x^2+2x+4)}}{4x^2} \right), \]
shown in figures \ref{fig:wrpm} and \ref{fig:r5},
which are tangent to the lines $y = \pm \frac{\sqrt{6}}{8} (x-2)$.

\begin{prop}
\label{prop:whomeo}
The functions $\ww_\pm$ are orientation preserving homeomorphisms
onto their respective images.
\end{prop}

\proof
The Jacobian matrix of $\tilde \ww_\pi$ is
\[ D\ww_\pm(x,y) = \begin{pmatrix}
\displaystyle\frac{2 (\omega_\pm)_x}{(1 - (\omega_\pm))^2}x +
\frac{1 + (\omega_\pm)}{1 - (\omega_\pm)} &
\displaystyle\frac{2 (\omega_\pm)_y}{(1 - (\omega_\pm))^2}x \\ \\
\displaystyle\pm\frac{(\omega_\pm)_x}{(1 + (\omega_\pm))^2}y &
\displaystyle\pm\frac{(\omega_\pm)_y}{(1 + (\omega_\pm))^2}y \pm
\frac{(\omega_\pm)}{1 + (\omega_\pm)} \end{pmatrix}, \]
with determinant given by
\[ \det D\ww_\pm(x,y) =
\pm\frac{1}{1 - \omega_\pm^2}
\left( \frac{2 (\omega_\pm)_x \omega_\pm x}{1 - \omega_\pm} +
(\omega_\pm)_y y + \omega_\pm (1 + \omega_\pm) \right). \]
It follows from lemmas \ref{lemma:omega} and \ref{lemma:signofomegax}
that each term in the sum between parenthesis have the same sign
and $\det D\ww_\pm(x,y) > 0$ if $y \ne 0$.

Points in the horizontal axis are fixed points of $\ww_\pm$.
Figure \ref{fig:wrpm} indicates
that the boundary of the domains are taken to simple closed curves,
which in turn implies the result, from standard degree theory.
A more rigorous and rather lengthy proof is possible
using estimates and a little topology, but will be omitted.
\qed

Notice that $\tilde \ww_\pi$ reverses orientations for $(x,y) \in R_-$
but $\ww_-$ preserves orientations.


\section{Quadratic convergence for  Wilkinson's iteration}

\begin{theo}
\label{theo:localantiwilk}
There exists an open neighborhood $\Aa \subset \ILa$ of $P_0$
and a closed set $\Xx \subset \Aa$ of zero measure,
invariant under $\ww$, on which the iteration converges quadratically to $P_0$.
The part of $\Xx$ with positive $y$ coordinate is homeomorphic
to the Cartesian product of a Cantor set and an open interval.
\end{theo}

\begin{figure}[ht]
\begin{center}
\psfrag{p0}{$p_0$}
\psfrag{30}{$(3,0)$}
\epsfig{height=35mm,file=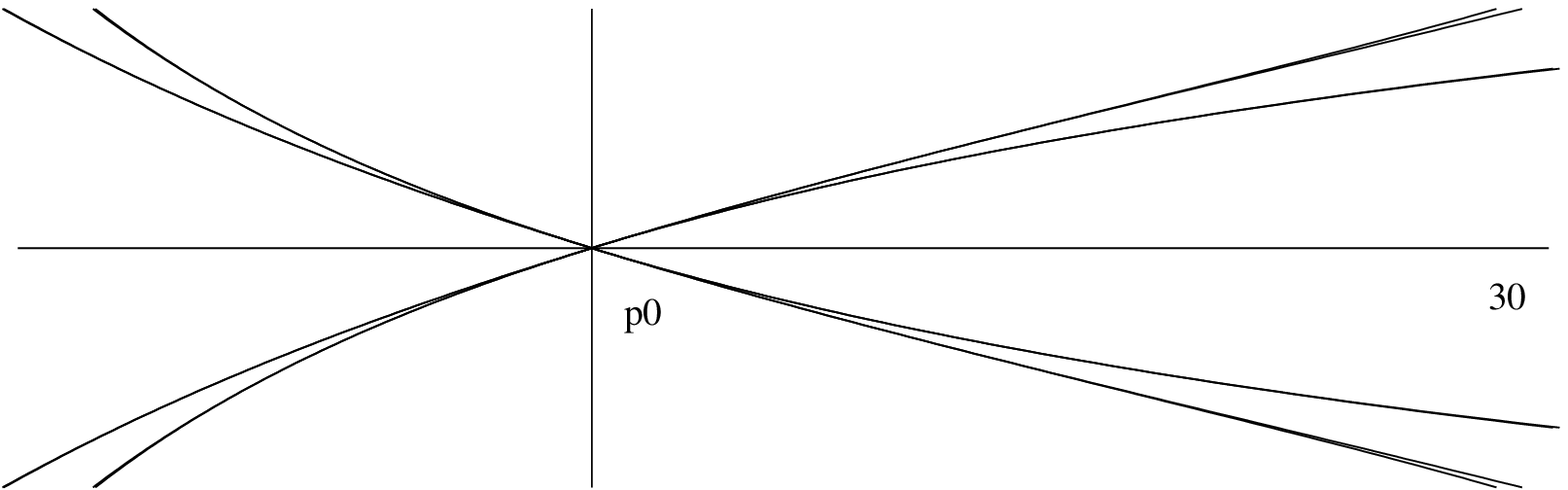}
\end{center}
\caption{\capsize The set $\Xx$ near $p_0$; in scale}
\label{fig:xxx}
\end{figure}

Numerical evidence indicates that we can take $\Aa = \UpLa$.
Figure \ref{fig:xxx} shows $\Xx$ in $\pi$-bidiagonal coordinates:
$\Xx$ and its mirror image at the $y$ axis map the set of {\it all}
matrices in $\ILa$ for which Wilkinson's iteration converges quadratically.
Since the Cantor set is extremely thin, the fine structure of the set $\Xx$
is invisible in the figure;
the curves $\ww_\pm^{-1}(r)$ fit inside the largest gaps of $\Xx$
in each quadrant.
As we shall see, an equivalent characterization of $\Xx$
is {\it wedge invariance}:
$\Xx$ is the set of points whose forward orbit under $\ww$
is eventually contained in $V_a$.
Propositions \ref{prop:antiwilk} and \ref{prop:Xquadratic}
below imply the theorem and provide additional,
more technical information about $\Xx$.

In a self-evident notation, we speak of the upper and lower half-wedges
and of the NE, NW, SE and SW faces of a wedge $V_a$.
Given $z_0 \in R$, set $z_{k+1} = \ww(z_k)$;
this is well defined unless $z_k \in r$.
The sequence $(z_k)_{k \in \NN}$ is the {\it $\ww$-orbit} of $z_0$.

\begin{lemma}
\label{lemma:wedgeinvariance}
For sufficiently small $a > 0$,
if $(x,y) \not\in V_a$, $|y| \le a$, then $\ww(x,y) \not\in V_a$.
Furthermore, a $\ww$-orbit tends to $p_0$ if and only if
it is eventually contained in a wedge.
\end{lemma}

\proof
Consider a short segment
in the upper half plane starting at $p_0$ with argument $\theta$.
It is easy to see that for $\theta > \pi -\arctan(\sqrt{6}/8)$,
the image under $\ww_+$ of this segment is a curve
tangent to the horizontal axis at $p_0$ and to
the left of the vertical line $r$.
Similarly, if $\theta < \pi -\arctan(\sqrt{6}/8)$,
the curve is to the right of $r$.
An example of this is $\ww_+(r)$, shown in figure \ref{fig:r5}.
We remind the reader that $-\sqrt{6}/8$ is the slope of $\ww_+^{-1}(r)$ at $p_0$.
Since $-1/6 < -\sqrt{6}/8 < 0$, the images of the NW face
and of the line $r$ are to the left and right, respectively, of the wedge.
A similar remark holds for $\ww_-$ and the NE face.

Near the horizontal axis, $|\omega|/(1+\omega)$ can be assumed to be smaller
than $1$ and therefore the absolute value of the second coordinate of $z_k$
is decreasing. The slope of the line joining $z_0 = (x,y)$ and $z_1$ is
\[ \frac{y}{\omega_\pm}
\frac{(1-\omega_\pm)(1 + \omega_\pm \mp \omega_\pm)}{-2x(1+\omega_\pm)}. \]
Since $|\omega_\pm| \le 2|y|$ (lemma \ref{lemma:signofomegax})
and the second fraction tends to $1/4$ when $(x,y)$ tends to $p_0$,
the slope can be assumed to have absolute value greater than $1/9$,
i.e., to be steeper than the faces of the wedge.
Thus, $V_a$ is further from $z_{k+1}$ than from $z_k$.
\qed

The set $\Xx$ can now be defined either
as the set of points whose orbit is eventually contained in a wedge or
as the set of points $z$ for which $\lim_{k \to \infty} \ww^k(z) = p_0$.

An {\it $L$-flat arc} in $V_a$ is the graph
$\Gamma \subset V_a$ of a $L$-Lipschitz function $f: I \to \RR$.

\begin{lemma}
\label{lemma:achata}
There exist a wedge $V_{a^\ast}$ and a positive constant $L^\ast < 1/6$
with the following properties.
Suppose $\Gamma^+_0$ is an $L^\ast$-flat arc in $V_{a^\ast}$,
with left endpoint belonging to the NW face of $V_{a^\ast}$
and right endpoint in the vertical line $r$.
Then $\ww_+(\Gamma^+_0)$ contains $\Gamma_1$,
also an $L^\ast$-flat arc in $V_{a^\ast}$
with left (resp. right) endpoint in the NW (resp. NE) face of $V_{a^\ast}$.
Moreover, such arcs are uniformly pushed towards the horizontal axis:
\[ \max_{(x,y) \in \Gamma_1} y < 1/4 \min_{(x,y) \in \Gamma^+_0} y. \]
Furthermore, $\ww_+$ stretches the horizontal coordinate.
More precisely, let the endpoints of $\Gamma_1$ be $\ww_+(x_\pm, y_\pm)$
and for $x \in [x_+, x_-]$, let $\phi(x)$ be the first coordinate
of $\ww_+(x,y)$ where $(x,y) \in \Gamma^+_0$; then $\phi'(x) > 1$ for all $x$.

An analogous statement holds for the action of $\ww_-$ on an $L^\ast$-flat
arc $\Gamma^-_0$ with endpoints now belonging to $r$
and the NE face of the wedge.

Symmetry with respect to the horizontal axis implies similar results
for arcs in the lower half wedge.
\end{lemma}

Notice that on smaller wedges, the lemma still holds for the same
Lipschitz constant but given a wedge, the Lipschitz constant cannot
be taken arbitrarily small.

\proof
We prove the statements concerning the action of $\ww_+$ on the upper half wedge,
the others being similar.

In order to control the slope of images of $L$-flat arcs,
we proceed to prove the following claim.
Given $L > 0$, there exists $a > 0$ such that if $(x,y) \in V_a$ then:
\begin{enumerate}
\item{the eigenvalues $\lambda_0$ and $\lambda_1$ of $D\ww_+(x,y)$
satisfy $|\lambda_0| < 1/4$, $\lambda_1 > 1/2$;}
\item{for the associated eigenvectors $v_i$,
$|\cot\arg v_0| < 1/L$ and $|\tan\arg v_1| < L$.}
\end{enumerate}
Indeed, from lemma \ref{lemma:signofomegax} and the formula for $D\ww_+$
in the proof of proposition \ref{prop:whomeo},
the entries in the second row of $D\ww_+$ tend to zero
when $(x,y)$ tends to $p_0$, the $(1,2)$ entry is bounded
and the $(1,1)$ entry is larger than $3/4$.
In a suggestive notation,
\[ D\ww_+ = \begin{pmatrix} a_1 & a_2 \\ \epsilon_1 & \epsilon_2 \end{pmatrix} \]
has eigenvalues
\[ \frac{(a_1 + \epsilon_2) \pm \sqrt{(a_1 - \epsilon_2)^2 + 4 a_2 \epsilon_1}}
{2}, \]
from which the estimates for $\lambda_0$ and $\lambda_1$ follow.
The eigenvectors can be written as $v_0 = C_0 (a_2, \lambda_0 - a_1)$
and $v_1 = C_1 (\epsilon_2 - \lambda_1, -\epsilon_1)$,
from which estimates for the arguments also follow,
completing the proof of the claim.

Assume without loss of generality that $L^\ast$ is so small that
\[ \max_{(x,y) \in \Gamma} y < 2 \min_{(x,y) \in \Gamma} y \]
for any $L^\ast$-flat arc $\Gamma$ in $V_{a^\ast}$.
Assume also that $\omega_+ < 1/8$ for all $(x,y) \in V_{a^\ast}$.
From
\[ \ww_+(x,y) = (x_1,y_1) = \left( \frac{1+\omega_+}{1-\omega_+} x,
\frac{\omega_+}{1+\omega_+} y \right) \]
we have $y_1 < y/8$ proving
\[ \max_{(x,y) \in \Gamma_1} y < 1/4 \min_{(x,y) \in \Gamma^+_0} y. \]
The claim implies that for all $(x,y) \in V_a$,
$v_1$ is in the east sector $|\arg v_1| < \arctan L$
and $v_0$ is in the north sector $\arctan L < \arg v_0 < \pi - \arctan L$.
Thus, the east sector is taken by $D\ww_+(x,y)$ to a subset of itself.
Setting $L = L^\ast$ and $a^\ast = a$,
this in turn implies that the image under $\ww_+$ of the arc $\Gamma^+_0$
in the statement of the lemma is an arc for which
the Lipschitz constant $L^\ast$ still holds.
As seen in the beginning of the proof, the endpoints of $\ww_+(\Gamma^+_0)$
are to the left and right of $V_{a^\ast}$.
The intersection of $\ww_+(\Gamma^+_0)$ with $V_{a^\ast}$ is $\Gamma_1$.

Write
\[ a_1 = \frac{1 + 2(\omega_\pm)_x x - \omega_\pm^2}{(1 - \omega_\pm)^2}. \]
Take $V_{a^\ast}$ so that $(1 - \omega_\pm)^2 < 1.01$,
$\omega_\pm^2 < 0.01$, $2(\omega_\pm)_x x > 0.1$ from which we learn
that $a_1 > 1.05$, completing the proof.
\qed

A {\it sign sequence} is a function $s: \NN \to \{\pm 1\}$
or $(s_0, s_1, s_2, \ldots)$.
The distance between two distinct sign sequences $s$ and $\tilde s$
is $3^{-k}$, where $k$ is the smallest number for which $s_k \ne \tilde s_k$.
There is a natural bi-Lipschitz homeomorphism between the set $\Ss$
of all sign sequences and the middle third Cantor set $\Kk_3 \subset [0,1]$:
take $s$ to $\sum_{k \ge 0} (1 + s_k)/3^{k+1}$.
More generally, a closed subset $\Kk$ of a Lipschitz graph $\Gamma$
is a {\it Cantor set} if it has empty interior (in the induced topology
in $\Gamma$) and no isolated points.
As is well known, a subset of a Lipschitz graph $\Gamma$ is a Cantor set
if and only if it is homeomorphic to $\Ss$.

Given $z_0 \in R$, the {\it $z_0$-sign sequence} $s^{z_0}$
is such that $s^{z_0}_k = +1$ iff $z_k \in R_+$ (where $z_{k+1} = \ww(z_k)$).

\begin{prop}
\label{prop:antiwilk}
Let $L^\ast$ and $a^\ast$ be as in lemma \ref{lemma:achata}.
Let $\Gamma_0$ be an $L^\ast$-flat arc in $V_{a^\ast}$
with endpoints in the NW and NE faces.
The set $\Xx \cap \Gamma_0$ is a Cantor set and the map taking
$z_0$ to the $z_0$-sign sequence is a bijection
from $\Xx \cap \Gamma_0$ to $\Ss$.

The set $\Xx$ is the disjoint union of graphs of Lipschitz
functions $f_s: [0, a^\ast] \to \RR$ (one for each sign sequence $s$)
taking $y_0$ to the $x$ coordinate of the unique point
in the intersection of $\Xx$ with the arc $y = y_0$ with sign sequence $s$.
\end{prop}

Numerical analysis gives, for example,
\[ f_{(+,+,+,+,\ldots)}(1/10) \approx
1.70831765759310579903646760761255776476753484977976.
\]

\begin{figure}[ht]
\begin{center}
\psfrag{a0}{$a_0$}
\psfrag{a1}{$a_1$}
\psfrag{a2}{$a_2$}
\psfrag{b0}{$b_0$}
\psfrag{b1}{$b_1$}
\psfrag{b2}{$b_2$}
\psfrag{gamma0p}{$\Gamma_0^+$}
\psfrag{gamma1m}{$\Gamma_1^-$}
\psfrag{gamma2m}{$\Gamma_2^-$}
\psfrag{r}{$r$}
\psfrag{p0}{$p_0$}
\psfrag{Wp}{$\ww_+$}
\psfrag{Wm}{$\ww_-$}
\epsfig{height=30mm,file=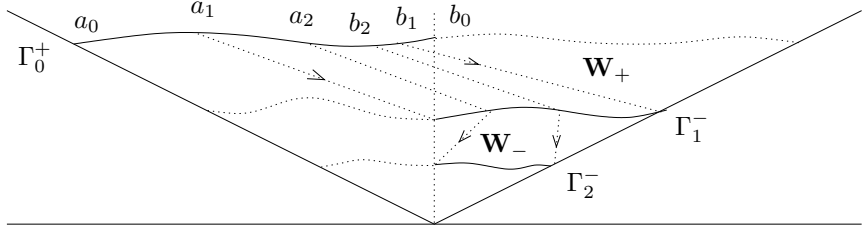}
\end{center}
\caption{\capsize Some curves $\Gamma_i$; schematic.}
\label{fig:gammagammagamma}
\end{figure}

\proof
Let $\Gamma_0^\pm = \Gamma_0 \cap R_\pm$.
From lemma \ref{lemma:achata}, the image of $\Gamma_0^\pm$ under $\ww_\pm$
contains an $L^\ast$-flat arc in $V_{a^\ast}$
with endpoints in the NW and NE faces.
For a sign sequence $s = (s_0, s_1, s_2, \ldots)$,
let $\Gamma_1$ be such an arc contained in $\ww_{s_0}(\Gamma_0^{s_0})$,
and, more generally, $\Gamma_{k+1}$ be an arc contained in
$\ww_{s_k}(\Gamma_k^{s_k})$.
Define intervals $I_k = [a_k, b_k] \subset \Gamma_0$
(see figure \ref{fig:gammagammagamma} for $s = (+,-,-,\ldots)$)
by \[ I_0 = \Gamma_0^{s_0}, \quad
\ww_{s_{k-1}} \cdots \ww_{s_1} \ww_{s_0} I_k = \ww^k I_k = \Gamma_k^{s_k}. \]

Again from lemma \ref{lemma:achata}, $|I_{k+1}| \le |I_k|/4$
and the intersection of the nested family of intervals $\cap_k I_k$
consists of the unique point in $\Xx \cap \Gamma_0$ with sign sequence $s$.
Thus the map from $\Ss$ to $\Xx \cap \Gamma_0$ taking $s$ to the point
with sign sequence $s$ is injective and continuous,
whence $\Xx \cap \Gamma_0$ is a Cantor set.

Now fix a sign sequence $s$.
Since the arc $y=y_0$ is $L^\ast$-flat, the function $f_s$ is well defined.
We show that $f_s$ is Lipschitz with constant $1/L^\ast$.
Indeed, assume by contradiction that $y_1$ and $y_2$ satisfy
\[ |f_s(y_1) - f_s(y_2)| > \frac{1}{L^\ast} |y_1 - y_2|; \]
the line through the points $(f_s(y_1), y_1)$ and $(f_s(y_2), y_2)$
has slope smaller than $L^\ast$ and is therefore $L^\ast$-flat.
Thus, there are two points on the intersection of $\Xx$ with
an $L^\ast$-flat arc with the same sign sequence, a contradiction.
\qed

The set $\Xx$ is rather thin, with Hausdorff dimension
(at least in a neighborhood of $p_0$) equal to $1$;
we do not present a proof of this fact.
Numerics suggests that $\Xx$ is a union of smooth curves $\Xx_s$, $s \in \Ss$,
parametrized by $(f_s(y),y)$.
The curve $\Xx_s$ is taken by $\ww$ to $\Xx_{s'}$,
where $s'$ is the left shift of $s$:
$s' = (s(1), s(2), s(3), \ldots)$.

\begin{prop}
\label{prop:Xquadratic}
For $z_0 = (x_0, y_0) \in \Xx$, there exists positive constants $c, C$
such that, for $z_n = (x_n, y_n)$,
$c |y_n|^2 \le |y_{n+1}| \le C |y_n|^2$,
i.e., the convergence of $z_n$ to $p_0$ is strictly quadratic.
On the other hand, for $z_0 \in R - \Xx$
the convergence of $z_n$ is strictly cubic.
\end{prop}

\proof
Recall that $y_{n+1} = \frac{|\omega|}{1+\omega} y_n$.
From lemma \ref{lemma:omega}, there exist positive constants
$c_1, C_1$ such that $c_1 |y| \le |\omega| \le C_1 |y|$
for any point $(x,y) \in V_a$.
Also, we may assume that $1/2 < 1+\omega < 2$ and therefore
\[ \frac{c_1}{2} |y|^2 \le \frac{|\omega|}{1+\omega} |y| \le 2C_1 |y|^2 \]
and the first claim follows.

If the limit point is some $p = (x, 0)$, $x \ne 0$,
then there exists positive constants $c_1, C_1$ such that,
in a neighborhood of $p$,
$c_1 |y|^2 \le |\omega| \le C_1 |y|^2$.
This follows from the fact that $\omega$ is smooth and even near $p$
we may write a Taylor expansion as in theorem \ref{theo:wilk}.
The second claim now follows easily.
\qed

\vfill\eject

\bigskip\bigskip\bigbreak

{

\parindent=0pt
\parskip=0pt
\obeylines

Ricardo S. Leite, Departamento de Matemática, UFES
Av. Fernando Ferrari, 514, Vitória, ES 29075-910, Brazil

\smallskip

Nicolau C. Saldanha and Carlos Tomei, Departamento de Matem\'atica, PUC-Rio
R. Marqu\^es de S. Vicente 225, Rio de Janeiro, RJ 22453-900, Brazil

\smallskip

rsleite@cce.ufes.br
nicolau@mat.puc-rio.br; http://www.mat.puc-rio.br/$\sim$nicolau/
tomei@mat.puc-rio.br

}

\end{document}